\documentclass[aos,MSNbibl,nameyear,dvips]{arximspdf}
\usepackage{multirow}
\usepackage{graphicx}

\doi{10.1214/13-AOS1175C} 
\referstodoi{10.1214/13-AOS1175}
\volume{42}
\issue{2}
\pubyear{2014}
\firstpage{483}
\lastpage{492}

\makeatletter
\newcommand{\AR}{\mathit{AR}}
\newcommand{\eqref}[1]{(\ref{#1})}
\newcommand{\overset}{\stackrel}
\def\argmin{\mathop{\arg\min}}
\def\var{\operatorname{var}}
\def\argmin{\mathop{\arg\min}}
\makeatother

\begin{document}
\begin{frontmatter}

\title{Discussion: ``A significance test for the lasso''}
\runtitle{Discussion}

\begin{aug}
\author{\fnms{Jianqing} \snm{Fan}\corref{}\ead[label=e1]{jqfan@princeton.edu}\thanksref{t1}}
\and
\author{\fnms{Zheng Tracy} \snm{Ke}\ead[label=e2]{zke@princeton.edu}\thanksref{t2}}
\runauthor{J. Fan and Z. T. Ke}
\affiliation{Princeton University}
\address{Department of Operations Research\\
\quad and Financial Engineering\\
Princeton University\\
Princeton, New Jersey 08544\\
USA\\
\printead{e1}\\
\phantom{E-mail:\ }\printead*{e2}} 
\pdftitle{Discussion of ``A significance test for the lasso''}
\end{aug}
\thankstext{t1}{Supported by NIH Grants R01-GM072611 and R01-GM100474
and NSF Grant DMS-12-06464.}
\thankstext{t2}{Supported by NIH Grant R01-GM072611 and NSF Grant
DMS-12-06464.}

\received{\smonth{12} \syear{2013}}



\end{frontmatter}

We wholeheartedly congratulate Lockhart, Taylor, Tibshrani and
Tibshrani on the stimulating paper, which provides insights into
statistical inference based on the lasso solution path. The authors
proposed novel covariance statistics for testing~the significance of
predictor variables as they enter the active set, which \mbox{formalizes} the
data-adaptive test based on the lasso path. The observation that
``shrinkage'' balances ``adaptivity'' to yield to an asymptotic
$\operatorname{Exp}(1)$ null distribution is inspiring, and the
mathematical analysis is delicate and intriguing.

Adopting the notation from the paper under discussion, the main results
are that the covariance statistics (Theorem~1)
%
\begin{equation}
\label{eq1} (T_{k_0+1}, T_{k_0+2}, \ldots, T_{k_0+d})
\overset{d} {\to} \bigl( \operatorname{Exp}(1), \operatorname{Exp}(1/2), \ldots,
\operatorname{Exp}(1/d) \bigr)
\end{equation}
for orthogonal designs, and under the global null model (Theorem~2),
$T_1 \overset{d}{\to} \operatorname{Exp}(1) $, and under the general
model (Theorem~3), $P(T_{k_0+1} \geq t) \leq\exp(-t) + o(1)$. These
remarkable results are derived under a number of critical assumptions
such as the normality, the sure screening [borrowing the terminology of
\citet{fan2008sure}] or model selection consistency of the lasso
path. As pointed out in \citet{FanLi}, lasso introduces biases
that are hard to account for. This together with the popularity of
lasso give rise to the importance of this work, which results in
informal statistical inference for the lasso. We welcome the
opportunity to make a few comments.

\section{Asymptotic null distributions}

\begin{figure}

\includegraphics{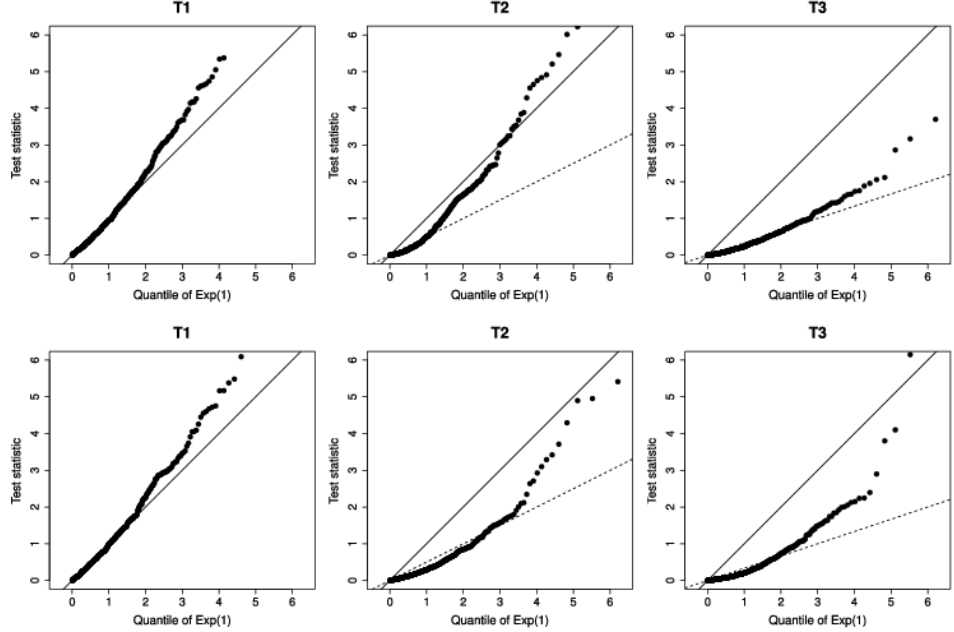}

\caption{Quantile--quantile plots of the covariance test statistics
versus their theoretical distributions under the global true null model
with ``equal correlation'' design (top panel) and ``$\AR(1)$'' design
(bottom panel) for
$n=100$ and $p=10$, based on 500 simulations.} \label{fig1}
\end{figure}

A natural question is how accurate the approximation \eqref{eq1} is
and whether it holds for more general design matrices. We illustrate
this using a small-scale numerical study. We take the same settings as
in Section~5.2 (Table~2) by considering the global null true model with
four types of design matrices: orthogonal, equal correlation, $\AR(1)$ and
block diagonal, where the parameter $\rho=0.8$. We fix $n=100$ and
$p=10$ and $50$. When $p= 50$, the marginal distributions of $\{T_1,
T_2, T_3\}$ are very close to the theoretical ones given by \eqref
{eq1}. However, when $p=10$, the approximation is not accurate for the
``equal correlation'' and ``$\AR(1)$'' designs. Figure~\ref{fig1} depicts
the results for $p=10$. The accuracies for the ``orthogonal'' and
``block diagonal'' designs are reasonable (omitted) and the accuracy
for $T_3$ is in general worse than those for $T_1$ and $T_2$.

%

\begin{figure}

\includegraphics{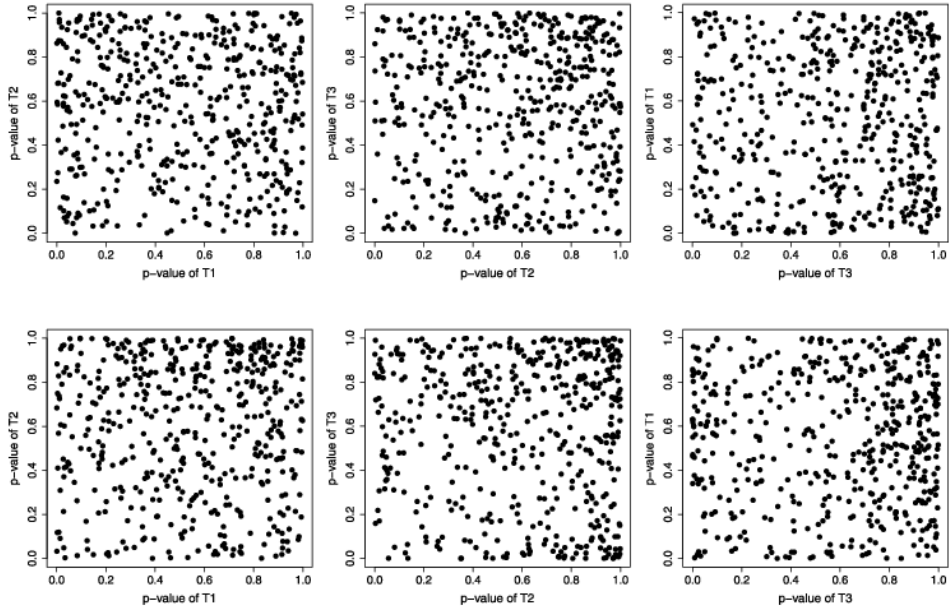}

\caption{$p$-values of $T_{k_0}$ versus those of $T_{k_0+1}$ for $n=100$
and $p=10$ with the orthogonal design matrices (top panel) and $n=100$
and $p=50$ with the equal-correlation design matrices (bottom panel)
based on $500$ simulations.} \label{fig2}
\end{figure}

To check the bivariate behavior of the covariance statistics $T_1$,
$T_2$ and $T_3$, we transform the statistics to have the asymptotic
uniform distribution using \eqref{eq1}. The scatter plots of those
transformed statistics are presented in Figure~\ref{fig2} based on 500
simulations. They are approximately uniformly distributed in the unit
square. This demonstrates that the test statistics are indeed
asymptotically independent and that the marginal distributions are
accurate for the given setting.

%

The simulation results presented in Figures~\ref{fig1} and~\ref{fig2} suggest that \eqref{eq1} holds for more general designs, not
just for orthogonal designs. This corresponds to suggesting that
Theorem~1 of the main paper holds more generally.

For a more general case in Theorem~3, the authors give a nice upper
bound. It requires a sure screening property and other conditions. A
large number of false positives in the set $A_0$ of the lasso path at
step $k_0$ should make the upper bound very crude and the upper bound
is tight when $A_0$ is model selection consistent.
This can easily be seen from the orthogonal design case with the global
null true model.
In this case, from \eqref{eq1},
%
\begin{equation}
T_{k_0+1} \stackrel{a} {\sim} \operatorname{Exp}(1)/k_0,
\end{equation}
which is of course stochastically bounded by $\operatorname{Exp}(1)$
but this bound can be very crude when $k_0$ is large.

Getting the sure screening property is difficult for lasso when the
irrespresentable condition [\citet{YuB}] does not hold. This was
demonstrated in \citet{fan2010sure} in which the design matrix is
generated such that $\{X_j\}_{j=1}^{p-50}$ are i.i.d. standard normal
variables and the last 50 predictor variables are
\[
X_k = \sum_{j=1}^s
\frac{(-1)^{j+1}}{5}X_j + \frac{\sqrt {25-s}}{5}\varepsilon_k,
\qquad k=p-49, \ldots, p,
\]
where $\{\varepsilon_k\}_{k=p-49}^p$ are i.i.d. standard normal
variables and $\{X_j\}_{j=1}^s$ are important variables. They also
noted that the larger the intrinsic model size $s$, the harder the
irrespresentable condition to hold; the larger the dimensionality, the
harder the condition. These follow from the definition of
irrespresentable condition. The question then arises what the null
distribution is when there are many false positives or even some false
negatives.

\begin{figure}

\includegraphics{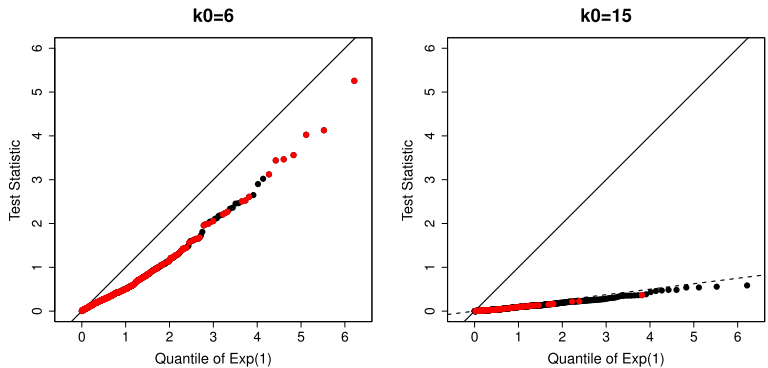}

\caption{Quantile--quantile plots of the covariance test statistics
versus $\operatorname{Exp}(1)$. By taking $k_0=6$ and $k_0=15$, the percentiles of sure
screening are 43.4\% and 87.8\%, respectively. The black/red dots
correspond to simulations with/without sure screening at step $k_0$.
The dash line on the right panel has slope $1/9$, matching the
distribution given by \protect\eqref{eq1} with $k_0=6$ and $d=9$.}
\label{fig3}
\end{figure}

To provide the insights, we fix $n=600$, $p=2000$, and $s=6$, take the
regression coefficient vector $\bolds{\beta}$ with $\beta
_1=\cdots=\beta_s=5$ and $\beta_{s+1}=\cdots=\beta_p=0$, and
simulated $500$ data sets. We computed the test statistic $T_{k_0+1}$
at $k_0=6$ and $k_0=15$. The results are shown in Figure~\ref{fig3}.
As expected, Theorem~3 continues to hold, but the bound is uselessly
crude. For $k_0 = 6$, there are only 43.4\% of the lasso paths that
have the sure screening or equivalently the model selection
consistency; others have both false positives and false negatives. As a
result, while Theorem~3 continues to hold, the bound is too crude. We
have also taken $k_0=15$, which makes 87.8\% of lasso paths to have
sure screening. In this case, there are many (at least 9) false
positives. Not knowing the true model size is 6, we compare it with
$\operatorname{Exp}(1)$ distribution, which shows again that Theorem~3
is correct, but the bound is too crude to be useful. Interestingly,
although this is not supported by Theorems 1--3, the test statistic
$T_{k_0+1}$ with $k_0 = 15$ is very close to $\operatorname
{Exp}(1/9)$, even though there are many false positives or even false
negatives. Is there any deeper theory underpinning the plot or is it
just a coincidence?

%

Another important condition is the normality assumption. This does not
seem as critical, thanks to the central limit theorem. For the
orthogonal design case, the variable $X_j^T y$ is approximately normal
under some mild conditions. For the logistic regression and Cox's
proportional hazards models, Figures~8 and 9 of the main paper show
that the covariance test statistic has approximately $\operatorname
{Exp}(1)$ distribution. Formal verifications of these results pose some
technical challenges, but are interesting research problems.


\section{Choice of the model size $k_0$}
The choice of model size $k_0$ is critically important. First, it
should be large enough to ensure the sure screening. Second, it should
not be too large to make overconservative inferences. For the current
paper, $k_0$ directly relates to the null distribution that is used for
computing $p$-values.

Let $T_k$ be the covariance statistic, defined by (7) and simplified in
(9) in the main paper. For a given $k_0$, define
\[
\label{eq4} \widetilde{T}_{k_0, j} = j T_{k_0+j}\qquad\mbox{for }
j = 1, \ldots, d.
\]
When $k_0$ is the correct model size so that the model selection
consistency holds, from (\ref{eq1}), $\{\widetilde{T}_{k_0, j}\}_{j=1}^d$
is a sequence of i.i.d. $\operatorname{Exp}(1)$ random variables. Therefore,
the average
%
\begin{equation}
\label{eq5} Q_{k_0} = \frac{1}{d} \sum
_{j=1}^d \widetilde{T}_{k_0, j} \approx1.
\end{equation}
%
A natural choice of $k_0$ is the one that makes $Q_{k_0}$
closest to its expected value $1$, namely
%
\begin{equation}
\label{eq5} \hat{k}_0 = \argmin_{k\ \mathrm{in\ a\ range}} |Q_k -
1|.
\end{equation}
The rationale is that when $k_0$ is the true model size, for example,
\[
E (Q_{k_0+1}) = \frac{1}{d} \sum_{j=1}^d
\frac{j}{j+1} = 1 - \frac
{1}{d}\sum_{j=1}^{d}(j+1)^{-1},
\]
which is less than $1$ and when $k < k_0$, $EQ_k$ is expected to be
much bigger than 1 (see Table~\ref{tab1}).

%
\begin{table}
\tabcolsep=0pt
\caption{Selection of the model size $k_0$. $n=500$ and the true
$k_0=2$. Based on $1000$ simulations, the mean of $Q_{k}$ (with
standard deviation in the parenthesis) and the distribution of selected
$\hat{k}_0$ are displayed} \label{tab1}
\begin{tabular*}{\tablewidth}{@{\extracolsep{\fill}}@{}lccccccc@{}}
\hline
$\bolds{p}$ & $\bolds{d}$ & $\bolds{k}$ & \textbf{0} & \textbf{1} & \textbf{2} & \textbf{3} & \textbf{4}\\
\hline
\phantom{00 }10 & \phantom{0}6 & mean $Q_{k}$ & 9.30 (2.3) & 4.40 (1.3) & 0.76 (0.43) & 0.48 (0.28) & 0.33 (0.22)\\
&& $\operatorname{prob} (\hat{k}_0=k$) & $0.0\%$ & $0.5\%$ & $79.9\%$ & $15.2\%$ & $4.4\%$
\\[3pt]
1000 & \phantom{0}6 & mean $Q_{k}$ & 6.31 (2.0) & 3.00 (1.1) & 0.93 (0.39) & 0.66 (0.29) & 0.53 (0.24)\\
&& $\operatorname{prob} (\hat{k}_0=k$) & $0.0\%$ & $0.04\%$ & $64.5\%$ & $20.5\%$ &$10.8\%$
\\[3pt]
1000 & 20 & mean $Q_{k}$ & 2.58 (0.62) & 1.53 (0.36) & 0.85 (0.20) & 0.72 (0.17) & 0.64 (0.16)\\
&& $\operatorname{prob} (\hat{k}_0=k$) & $0.0\%$ & $22.3\%$ & $61.6\%$ & $11.6\%$ & $4.5\%$\\
\hline
\end{tabular*}
\end{table}

To see the accuracy of this method, we note that it is typically the
hardest to differentiate the choice of $k_0$ and $k_0+1$ when the true
model size is $k_0$. The variance of the difference is
\begin{eqnarray*}
\var(Q_{k_0} - Q_{k_0+1}) & = & d^{-2}
\var(T_{k_0+1} + \cdots+ T_{k_0+d} - d T_{k_0+d+1})
\\
& = & d^{-2} \bigl(1+2^{-2} + \cdots+ d^{-2} +
d^2/(d+1)^2\bigr).
\end{eqnarray*}
%
It follows that
\begin{eqnarray*}
\frac{E(Q_{k_0}-Q_{k_0+1})}{\var(Q_{k_0} - Q_{k_0+1})^{1/2}}&=& \frac
{\sum_{j=1}^{d} (j+1)^{-1}}{(1+2^{-2} + \cdots +d^{-2} +
d^2/(d+1)^2)^{1/2}}
\\
& \asymp&\frac{\log(d)}{\sqrt{1+\pi^2/6}}\qquad\mbox{as
} \,d\to \infty.
\end{eqnarray*}
Thus, the signal to noise ratio is large when $d$ is large, but
increases slowly with~$d$. Therefore, in practice, we do not wish to
take a too large $d$ due to the accuracy of approximation \eqref{eq1}.

We conducted a numerical experiment where $n=500$ and $p=10$ and
$1000$. The predictors $\{X_j\}_{j=1}^p$ are i.i.d. standard normal
variables. Let $\beta=(6, 6, 0, \ldots, 0)^T$, so the true $k_0=2$.
For fixed $d=6$ and $d=20$ (only when $p=1000$), we selected $\hat
{k}_0$ from $\{0, \ldots, 4\}$ to minimize $|Q_{k}-1|$.
Table~\ref{tab1} summarizes the results based on $1000$ simulations.
When $p=10$, the percentiles of $\hat{k}_0=2$ (accurate) and $\hat
{k}_0=3$ (overshooting by $1$) are about $80\%$ and $15\%$, and there
are almost no undershootings ($\hat{k}_0<2$). When $p=1000$, the
accuracy decreases to about $65\%$, but there are still almost no
undershootings. Interestingly, when we increase $d$ to $20$, the
results become inferior, with about $22\%$ of undershootings. This
suggests that $d$ should not be chosen too large that smooths out the
signals in $Q_k$ for $k < k_0$ and makes \eqref{eq1} inaccurate.


\section{Power of the tests}

When we test the first few variables that enter the active set of
lasso, it is very often that there remain true active variables not yet
selected. The covariance test statistics are just one of many
possibilities, constructed carefully and intriguingly and supported by
the nice asymptotic null distribution. However, they are not
necessarily the most powerful tests.

To understand the possible loss of the power of the covariance test, we
consider again the simplest setting where the design matrix is
orthogonal, $k_0=0$ (so the null hypothesis is the global null) and
$\sigma=1$. It follows that
\[
T_1 = V_1 (V_1 - V_2),
\]
where $V_1$ and $V_2$ are the first and second largest elements of $\{
|X_j^Ty|\dvtx  1\leq j\leq p\}$.
The factor $V_1 - V_2$ makes the null distribution very beautiful, but
this can also reduce the power of the statistic $V_1$, which is equal
to square root of the maximum drop in RSS.

To demonstrate this, consider the specific alternative
\[
\beta_1=\beta_2=\theta, \qquad\beta_3=
\cdots=\beta_p=0,
\]
where $\theta\gg\sqrt{\log(p)}$. With probability tending to $1$,
$|X_1^Ty|$ and $|X_2^Ty|$ are the largest two elements. As a result,
$V_1$ is stochastically equivalent to that of $(\theta+\max\{\epsilon
_1,\epsilon_2\}) + o_p(1)$ and $V_1 - V_2 = |\varepsilon_1 -
\varepsilon_2| + o_p(1)$ with $\epsilon_1,\epsilon_2$ being
independent standard normal variables. It follows that
%
\begin{equation}
\label{eq3} T_1/\theta\overset{d} {\to} \bigl|N(0,\sqrt{2})\bigr|\quad\mbox{and}\quad V_1/ \theta\overset{d} {\to} 1.
\end{equation}
The statistic $T_1$ and the maximum drop of RSS $V_1^2$ indeed have
asymptotic power one. On the other hand, \eqref{eq3} shows clearly
that $T_1$ is corrupted by an extra noise $|N(0, \sqrt{2})|$ and is
therefore less powerful.

%
\begin{figure}

\includegraphics{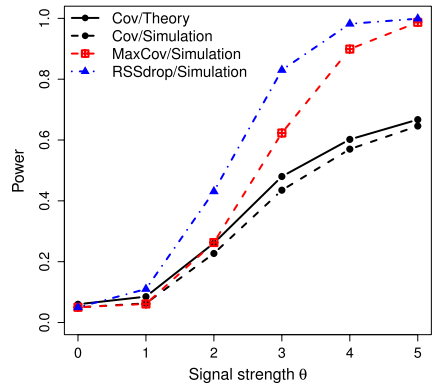}

\caption{Power curves based on $1000$ simulations. $n=100$, $p=10$ and
the predictors are drawn i.i.d. from $N(0,1)$. ``Cov/Theory'' and
``Cov/Simulation'' refer to the covariance test statistic $T_1$, with
critical value being the $95\%$ quantile of\, $\operatorname{Exp}(1)$
and the sample $95\%$ quantile. ``MaxCov/Simulation'' refers to the
maximum of $T_1$ and $T_2$, and ``RSSdrop/Simulation'' refers to the
maximum drop in RSS, with critical values being the sample $95\%$ quantile.}
\label{fig4}
\end{figure}

We illustrate this point using a small-scale numerical study. We use
similar settings as the left panel of Figure~4 in the main paper
($n=100$, $p=10$ and ``orthogonal design''). Instead of having only one
truly nonzero regression coefficient, we set two equal nonzero
regression coefficients. Figure~\ref{fig4} shows the estimated power
curves. When there is only one true nonzero coefficient, the covariance
test statistic and the maximum drop in RSS have similar powers as shown
in Figure~4 of the main paper. On the other hand, when there are two
equal nonzero coefficients, the statistic of maximum drop in RSS has a
larger power, especially when the signal strength $\theta$ is large.
Interestingly, when we compute the covariance test statistics in this
case, it is quite often that the first entering variable is not very
significant but the second one is. We also compute the power when
looking at the maximum of $T_1$ and $T_2$. It turns out that this test
is more powerful than using $T_1$ only. See Figure~\ref{fig4}.

\section{Validity of the results to other penalties}

A natural question is whether or not the results in the paper are tied
to the lasso path. Given many nice bias properties of folded concave
penalty [\citet{FanLi}] and weighted lasso penalty [\citet
{Zou}] functions, it is natural to examine the solution paths created
by those penalty functions.

For a general penalty function $p_{\lambda}(\cdot)$, we define the
covariance test statistic at the knot $\lambda_{k}$ the same as (5) in
the main paper,
except that $\hat{\beta}(\lambda_{k+1})$ and $\tilde{\beta
}_A(\lambda_{k+1})$ are computed with $\Vert\beta\Vert_1$ replaced
by $\sum_j p_{\lambda}(\beta_j)$ in the expressions. Although there
are the issues on the uniqueness of the folded concave penalized
least-squares, \citet{fan2011nonconcave} show that folded concave
penalized least-squares estimator is indeed unique in the sense of
restricted global optimality.

As in the main paper, we examine the showcase example in which the
design matrix is orthogonal. In this case, the penalized least-squares
with folded concave penalty is unique [\citet{FanLi}]. By direct
calculation,
%
\begin{equation}
\label{equ5} T_k = V_k\cdot h_{V_{k+1}}(V_k)/
\sigma^2,
\end{equation}
where $h_{\lambda}(\cdot)$ is a thresholding function defined by
$h_{\lambda}(x)= \arg\min_u\{\frac{1}{2}(u-x)^2 + p_{\lambda}(u)\}
$. For the SCAD penalty [\citet{FanLi}] with a parameter $a>2$,
\[
T_k^{\mathrm{scad}} = \cases{
\displaystyle V_k(V_{k}-V_{k+1})/\sigma^2, &\quad $V_k\leq2V_{k+1}$,
\vspace*{5pt}\cr
\displaystyle\frac{a-1}{a-2} V_k \biggl(V_{k}-\frac{a}{a-2}V_{k+1} \biggr)\Bigm/\sigma ^2, &\quad $2V_{k+1} < V_k <aV_{k+1}$,
\vspace*{5pt}\cr
\displaystyle V_k^2/\sigma^2, &\quad $V_k\geq a V_{k+1}$.}
\]
%

We can similarly show that for any fixed $k\geq1$,
%
\begin{equation}
\label{eq6} \bigl(T_1^{\mathrm{scad}}, T_2^{\mathrm{scad}},\ldots, T_k^{\mathrm{scad}} \bigr) \overset{d} {\to } \bigl(
\operatorname{Exp}(1),\operatorname{Exp}(1/2), \ldots, \operatorname{Exp}(1/k)
\bigr),
\end{equation}
under the global null true model.
\begin{pf}
Let $F(x)=(2\Phi(x)-1)I\{x>0\}$. From the proof of Lemma~3 in the main
paper, for $a_p=F^{-1}(1-1/p)$ and $b_p=pF'(a_p)$, the random variables
$W_0=b_p(V_{k+1}-a_p)$ and $W_i=b_p(V_i-V_{i+1})$, $i=1,\ldots,k$,
converge jointly:
%
\begin{equation}
\label{eq7} (W_0, W_1, W_2, \ldots,
W_k) \overset{d} {\to} ( -\log G_0, E_1,
E_2/2,\ldots, E_k/k ),
\end{equation}
where $G_0, E_1, \ldots, E_k$ are independent, $G_0$ is Gamma
distributed with scale parameter $1$ and shape parameter $k$, and
$E_1,\ldots, E_k$ are standard exponentials. In addition, $a_p, b_p\to
\infty$ and $a_p/b_p\to1$ as $p\to\infty$.

Note that $T_i^{\mathrm{scad}}=T_i^{\mathrm{lasso}}$, $i=1,\ldots, k$, on the event
$B=\{ V_i\leq2V_{i+1}, 1\leq i\leq k\}$. By \eqref{eq7} and the fact
that $a_pb_p\to\infty$,
\[
P(V_i>2V_{i+1}) = P \Biggl( iE_i+\log
G_0-\sum_{j=i+1}^k
jE_j > a_pb_p \Biggr) = o(1).
\]
Then $P(B^c)\leq\sum_{i=1}^k P(V_i>2V_{i+1})=o(1)$. Therefore, \eqref
{eq6} follows immediately from Lemma~3 and the Slucky's lemma.
\end{pf}


For the MCP penalty [\citet{zhang2010nearly}] with a parameter
$\gamma>1$, it can be shown similarly that
\[
\bigl(T_1^{\mathrm{mcp}}, T_2^{\mathrm{mcp}}, \ldots,
T_k^{\mathrm{mcp}} \bigr) \overset{d} {\to} \frac{\gamma}{\gamma-1} \bigl(
\operatorname{Exp}(1),\operatorname {Exp}(1/2), \ldots, \operatorname{Exp}(1/k)
\bigr).
\]
For the weighted lasso penalty [\citet{Zou}],
the solution path depends on order statistics of variables $\{
w_j^{-1}|X_j^Ty|\}_{j=1}^p$, where $w_j$ is the weight for variable
$j$. These variables are not identically distributed. It remains an
interesting question to what extent the current results can be generalized.

\section{Further comments}

The mathematical results are derived when $d$ and $k_0$ are finite. A
more interesting asymptotic framework is to let both $d$ and $k_0$
diverge with $n$. See, for example, \citet{fan2011nonconcave} for
the joint asymptotic distribution when the dimensionality grows with
sample size.

The beautiful results in the paper are derived under the assumptions
that the signals are very strong and the designs are so nice that sure
screening is possible. These assumptions are difficult to meet in
practice. Even when they are met, we need to specify $k_0$ which is
hoped to be small and contains all important variables (sure
screening). Sure screening assumption implies that the null hypothesis
is true. What are we testing: sure screening hypothesis or significance
of the newly entered variable? Under the sure screening assumption, why
not run the least-squares based on the screened predictors and use
splitted data (when needed), as suggested in \citet{fan2008sure}
and \citet{Wasserman}? The statistical inference can be based
upon the low-dimensional least-squares theory. To utilize the
asymptotic null distribution without conservatism, we need to have the
model selection consistency assumption: the first $k_0$ variables
contain all important variables. If so, why do we need the significance
tests of the newly entered variables? Relaxing the model selection
consistency to sure screening does not help the matter very much. Using
the standard exponential distribution as the upper bound of the
$p$-values, we can mislabel many ``important variables'' as ``unimportant
ones,'' a missed discovery that we strive to avoid in high-dimensional inference.

The authors mentioned in the paper that they plan to construct
confidence regions for the lasso $\hat{\bolds{\beta}}(\lambda
)$ at specific $\lambda$. The challenge here is that there are biases
involved in the lasso fit. Another challenge is to give a formal
confidence assessment that a group of ``unimportant variables'' are
really unimportant. The efforts are certainly welcome [see, e.g.,
\citet{meinshausen2009p} and \citet{ZhangZhang11}]. We
would like to note that for the folded-concave penalized least-squares
or likelihood, the resulting estimator is the oracle estimator with
probability tending to one [\citet{fan2011nonconcave}].
Therefore, the confidence intervals can easily be constructed based on
the low-dimensional likelihood inference. However, it also remains to
give confidence assessment that a group of ``unimportant variables''
are really unimportant.

The authors have mentioned a couple of times the null distribution of
the largest RSS drop. This is equivalent to $(1 - \gamma_n^2)$, where
$\gamma_n$ is the maximum correlation coefficient between the
residuals at the current step of the forward regression and the
covariates. Under the global null true model, this is the maximum
spurious correlation between the response and each variable. The
asymptotic distribution for the maximum spurious correlation in case
where all predictors are independent has been derived in \citet
{cai2011limiting}. However, we do not expect that the asymptotic null
distribution is accurate enough for many applications.

In conclusion, the idea and results in the main paper are insightful
and amazing. The technical arguments are delicate and ingenious. The
authors should be congratulated again for successful adaptive inference
based on the lasso solution path. We hope that our comments contribute
positively to the understanding of this seminal article.



%

\printaddresses

\end{document}